\documentclass[12pt,oneside,english]{amsart}
\usepackage{newcent}
\usepackage[T1]{fontenc}
\usepackage[latin1]{inputenc}
\usepackage{amssymb}

\makeatletter
 \theoremstyle{plain}    
 \newtheorem{thm}{Theorem}[section]
 \numberwithin{equation}{section} 
 \numberwithin{figure}{section} 
 \theoremstyle{plain}
 \theoremstyle{definition}
 \newtheorem{defn}[thm]{Definition}
 \theoremstyle{plain}    
 \newtheorem{prop}[thm]{Proposition} 
 \theoremstyle{definition}
  \newtheorem{example}[thm]{Example}
 \theoremstyle{plain}    
 \newtheorem{lem}[thm]{Lemma} 

\address{Department of Mathematics, UCLA,
Los Angeles, CA 90095, USA}
\email{shlyakht@math.ucla.edu}

\AtBeginDocument{
  
}

\usepackage{babel}
\makeatother
\begin{document}

\title{Notes on Free Probability Theory}

\author{Dimitri Shlyakhtenko.}

\maketitle
These notes are from a 4-lecture mini-course taught by the author
at the conference on von Neumann algebras as part of the {}``Géométrie
non commutative en mathématiques et physique'' month at CIRM. 

\tableofcontents{}

\section{Free Independence and Free Harmonic Analysis.}

Free probability theory was developed by Voiculescu as a way to deal
with von Neumann algebras of free groups. In addition to the view
of von Neumann algebras as {}``non-commutative measure spaces'',
which was already presented in this conference, free probability theory
considers von Neumann algebras as {}``non-commutative probability
spaces''. 

There are by now several standard references on free probability theory,
of which we mention two: \cite{DVV:book,dvv:lectures}.

\subsection{Probability spaces.}

Recall that a \emph{classical probability space} is a measure space
$(X,B,\mu)$. Here $B$ is a sigma-algebra of subsets of $X$, and
$\mu$ is a measure, which is a \emph{probability measure}, i.e. $\mu(X)=1$.
One thinks of $X$ as a set of events and for $Y\in B$, the measure
$\mu(Y)$ is a probability of an event occurring in the set $Y$.

\subsubsection{Random variables; laws.}

An alternative point of view on probability theory involves considering
\emph{random variables}, i.e., measurable functions $f:X\to\mathbb{C}$.
One can think of a random variable as a measurement, which assigns
to each event $x\in X$ a value $f(x)$. Note that the probability
of the value of $f$ lying in a set $A\subset\mathbb{C}$ is exactly
$\mu(f^{-1}(A))=(f_{*}\mu)(A)$. Thus the \emph{law}~of $f$, $\mu_{f}$,
defined to be the push-forward measure $\mu_{f}=f_{*}\mu$ on $\mathbb{C}$,
measures the probabilities that $f$ assumes various values.

\subsubsection{The expectation $E$.}

Let us say that $f\in L^{\infty}(X,\mu)$ is an essentially bounded
random variable. Then the integral\[
E(f)=\int f(x)d\mu(x)\]
has the meaning of the expected value of $f$. For this reason, the
linear functional $E:L^{\infty}(X,\mu)\to\mathbb{C}$ given by integration
against $\mu$ is called an expectation. We note that $E$ satisfies:
$E(1)=1$ (normalization), $E(f)\geq0$ if $f\geq0$ (positivity).

Note that the knowledge of $(X,B,\mu)$ is equivalent (up to an isomorphism
and up to null sets) to the knowledge of $L^{\infty}(X,\mu)$ and
$E$. Thus the notion of a classical probability space can be phrased
entirely in terms of commutative (von Neumann) algebras.

\subsection{Non-commutative probability spaces. }

We now play the usual game of dropping the word {}``commutative''
in a definition:

\begin{defn}
An algebraic non-commutative probability space is a pair $(A,\phi)$
consisting of a unital algebra $A$ and a linear functional $\phi:A\to\mathbb{C}$,
so that $\phi(1)=1$.
\end{defn}
Thus we think of $a\in A$ as a {}``non-commutative random variable'',
$\phi(a)$ as its {}``expected value'' and so on. Of course, any
classical probability space is also a non-commutative probability
space. But there are many interesting genuinely non-commutative probability
spaces. For example, if $\Gamma$ is a discrete group, we could set
$A=\mathbb{C}\Gamma$ (the group algebra) and $\phi=\tau_{\Gamma}$
(the group trace). Here if $g\in\Gamma\subset\mathbb{C}\Gamma$, then
$\phi(g)=0$ if $g\neq1$ and $\phi(1)=1$. The same construction
works with $A$ replaced by the reduced group $C^{*}$-algebra of
$\Gamma$, or the von Neumann algebra of $\Gamma$.

\subsubsection{Positivity.}

Operator algebras give one a {}``test'' of which algebraic non-commutative
spaces {}``exist in nature''. These are precisely those non-commutative
probability spaces that can be represented by (possibly unbounded)
operators on a Hilbert space $H$, so that $\phi$ is a linear functional
given by a vector-state, $\phi(a)=\langle h,ah\rangle$ for some $h\in H$.
If $A$ is a $*$-algebra, it is not hard to characterize these (via
the GNS construction) in terms of the properties of $\phi$: $\phi$
must be \emph{positive}, i.e., $\phi(a^{*}a)\geq0$ for all $a\in A$.

\subsubsection{The law of a random variable.}

Recall that we assigned to a classical random variable $f$ its law
$\mu_{f}$. If $A$ is an algebra of operators on a Hilbert space
$H$, $\phi(\cdot)=\langle h,\cdot h\rangle$ and $a\in A$ is self-adjoint,
then the spectral theorem gives us a measure $\nu_{a}$ on $\mathbb{R}$
valued in the set of projections on $H$, so that\[
a=\int td\nu_{a}(t).\]
If we let\[
\mu_{a}=\phi\circ\nu_{a},\]
then $\mu_{a}$ is a measure on $\mathbb{R}$. It is not hard to check
that if we are in the classical situation and $a\in L^{\infty}(X,\mu)$,
$H=L^{2}(X,\mu)$, $h=1$, then this construction gives us precisely
the law of $a$.

\subsubsection{Moments.}

However, if $a$ is not self-adjoint, or if we are dealing with a
$k$-tuple of random variables, there is no description of the law
of $a$ in terms of a measure.

Fortunately, for $f\in L^{\infty}(X,\mu)$ the \emph{moments}~of
$f$, i.e., the expected values $E(f^{p})$, $p=1,2,\ldots$ are exactly
the same as the moments of the law $\mu_{f}$ of $f$. Indeed,\[
E(f^{p})=\int t^{p}d\mu_{f}(t)\]
is exactly the $p$-th moment of $\mu_{f}$. For essentially bounded
$f$, the moments of $\mu_{f}$ determine $\mu_{f}$.

Thus given a family $F$ of variables $a_{1},\ldots,a_{n}\in A$,
we say that an expression of the form $\phi(a_{i_{1}}\cdots a_{i_{p}})$
is the \emph{$i_{1},\ldots,i_{p}$-th moment} of the family $F$.
The collection of all moments can be thought of as a linear functional
$\mu_{F}$ defined on the algebra of polynomials in $n$ indeterminates
$t_{1},\ldots,t_{n}$ by\[
\mu_{F}(p)=\phi(p(a_{1},\ldots,a_{n})).\]
This functional $\mu_{F}$ is called the \emph{joint law}, \emph{}or
\emph{joint distribution}, of the family $F$.

\subsection{Classical independence.}

\begin{defn}
Two random variables $f$ and $g$ in $L^{\infty}(X,\mu)$ are called
independent, if\[
E(f^{n}g^{m})=E(f^{n})E(g^{m})\]
for all $n,m\geq0$.

Equivalently, $E(FG)=0$ whenever $E(F)=E(G)=0$ and $F$ is in the
algebra $W^{*}(f)$ generated by $f$, while $G\in W^{*}(g)$.
\end{defn}
The equality $E(fg)=E(f)E(g)$ is a consequence of the statement that
{}``the probability that the value of $f$ lies in a set $A$ and
the value of $g$ lies in the set $B$ is the product of the probabilities
that the value of $f$ lies in $A$ and the value of $g$ lies in
$B$'', which is a more familiar way of phrasing independence.

If $X=X_{1}\times X_{2}$ and $\mu=\mu_{1}\times\mu_{2}$, then any
functions $f,g$ so that $f$ depends only on the $X_{1}$ coordinate
and $g$ only on the $X_{2}$ coordinate are independent. Note that
another way of saying this is that the random variables $f_{1}\otimes1$
and $1\otimes g_{1}$ in $L^{\infty}(X_{1},\mu_{1})\bar{\otimes}L^{\infty}(X_{2},\mu_{2})$
are independent, for any $f_{1}\in L^{\infty}(X_{1},\mu_{1})$ and
$g_{1}\in L^{\infty}(X_{2},\mu_{2})$. Thus independence has to do
with the operation of taking tensor products of probability spaces.

\subsection{Free products of non-commutative probability spaces.}

There is {}``more room'' in the non-commutative universe to accommodate
a different way of combining two non-commutative probability spaces:
free products. Just like the notion of a tensor product can be used
to recover the notion of independence, free products have led Voiculescu
to discover the notion of \emph{free independence}.

\subsubsection{Free products of groups.}

We start with a motivating example. Let $\Gamma_{1}$ and $\Gamma_{2}$
be two discrete groups. View the group algebra of the free product
$\mathbb{C}(\Gamma_{1}*\Gamma_{2})$ as a non-commutative probability
space by letting $\phi$ be the group trace; for $g\in\Gamma_{1}*\Gamma_{2}$,
$\phi(g)=0$ unless $g=1$.

Let us understand the relative positions of $\mathbb{C}\Gamma_{1}$
and $\mathbb{C}\Gamma_{2}$ inside of the group algebra of the free
product $\mathbb{C}(\Gamma_{1}*\Gamma_{2})$. Let $w\in\Gamma_{1}*\Gamma_{2}$
be a word. Thus $g=g_{1}\cdots g_{n}$ with $g_{j}\in\Gamma_{i(j)}$.
We may carry out multiplications and cancellations until we reduce
the word so that consequent letters lie in different groups; i.e.,
$i(1)\neq i(2)$, $i(2)\neq i(3)$ and so on. The resulting word is
non-trivial if all $g_{1},\ldots,g_{n}$ are non-trivial. Thus:\[
\phi(g_{1}\cdots g_{n})=0\]
provided that $g_{j}\in\Gamma_{i(j)}$, $i(1)\neq i(2)$, $i(2)\neq i(3)$,
$\ldots$, and $\phi(g_{1})=\phi(g_{2})=\ldots=0$.

By linearity we get: 

\begin{prop}
if $a\in\mathbb{C}(\Gamma_{1}*\Gamma_{2})$ has the form\[
a=a_{1}\cdots a_{n},\]
with $a_{j}\in\mathbb{C}\Gamma_{i(j)}$, $i(1)\neq i(2)$, $i(2)\neq i(3)$,
$\ldots$, and $\phi(a_{1})=\phi(a_{2})=\ldots=0$, then\[
\phi(a)=0.\]

\end{prop}
We now note that this proposition allows one to compute $\phi$ on
$\mathbb{C}(\Gamma_{1}*\Gamma_{2})=\mathbb{C}\Gamma_{1}*\mathbb{C}\Gamma_{2}$
in terms of its restriction to $\mathbb{C}\Gamma_{1}$ and $\mathbb{C}\Gamma_{2}$.
Indeed, an arbitrary element of $\mathbb{C}\Gamma_{1}*\mathbb{C}\Gamma_{2}$
is a linear combination of $1$ and of terms of the form\[
a_{1}\cdots a_{n},\qquad a_{j}\in\mathbb{C}\Gamma_{i(j)},\quad i(1)\neq i(2),\  i(2)\neq i(3),\ldots.\]
But then the equation\[
0=\phi((a_{1}-\phi(a_{1}))(a_{2}-\phi(a_{2}))\cdots(a_{n}-\phi(a_{n})))\]
allows one to express $\phi(a_{1}\cdots a_{n})$ in terms of values
of $\phi$ on shorter words. By induction, this allows one to express
$\phi$ in terms of $\phi|_{\mathbb{C}\Gamma_{1}}$ and $\phi|_{\mathbb{C}\Gamma_{2}}$.

\subsubsection{Free products of algebras.}

Such an expression is universal and works in any free product of two
algebras (not necessarily of group algebras). We thus say:

\begin{defn}
\cite{DVV:free} Let $(A_{1},\phi_{1})$ and $(A_{2},\phi_{2})$ be
two non-commutative probability spaces. We call the unique linear
functional $\phi$ on $A_{1}*A_{2}$ which satisfies\begin{eqnarray*}
\phi(a_{1}\cdots a_{n}) & = & 0,\quad a_{j}\in A_{i(j)},\quad i(1)\neq i(2),\  i(2)\neq i(3),\ldots,\\
 &  & \qquad\phi_{i(j)}(a_{j})=0,\quad\forall j\end{eqnarray*}
the free product of $\phi_{1}$ and $\phi_{2}$. It is denoted $\phi_{1}*\phi_{2}$.
\end{defn}
One can check that the free product of two positive linear functionals
is positive (to do so it is the easiest to make sense of the product
of the underlying GNS representations). Thus one can talk about (reduced)
free products of $C^{*}$-algebras or von Neumann algebras by passing
to the appropriate closure in the GNS representation associated to
the free product functional.

\subsubsection{Free independence.}

By analogy with the relationship between classical independence and
tensor products, Voiculescu gave the following definition:

\begin{defn}
\cite{DVV:free} Let $F_{1},F_{2}\subset(A,\phi)$ be two families
of non-commutative random variables. We say that $F_{1}$ and $F_{2}$
are \emph{freely independent}, if\[
\phi(a_{1}\cdots a_{n})=0\]
whenever $a_{j}\in\textrm{Alg}(1,F_{i(j)})$, $i(1)\neq i(2)$, $i(2)\neq i(3),$
$\ldots$, and $\phi(a_{1})=\phi(a_{2})=\ldots=0$.
\end{defn}
Here $\textrm{Alg}(S)$ denote the algebra generated by a set $S$.

We should point out a certain similarity between this definition and
the classical independence, where the requirement was that $E(FG)=0$
if $E(F)=E(G)=0$.

\subsection{Free Fock space.}

We give an example of freely independent random variables that does
not come from groups.

\subsubsection{Free Fock space.}

Let $H$ be a Hilbert space, $\Omega$ be a vector, and let\[
F(H)=\mathbb{C}\Omega\oplus H\oplus H\otimes H\oplus\cdots\]
be the Hilbert space direct sum of the tensor powers of $H$ (the
one-dimensional space $\mathbb{C}\Omega$ is thought of as the zeroth
tensor power of $H$). This space is called the free (or full) Fock
space, by analogy with the symmetric and anti-symmetric Fock spaces
(where the symmetric or anti-symmetric tensor product is used instead).

\subsubsection{Free creation operators.}

For $h\in H$ consider the left creation operator\[
\ell(h):F(H)\to F(H)\]
given by\[
\ell(h)h_{1}\otimes\cdots\otimes h_{n}=h\otimes h_{1}\otimes\cdots\otimes h_{n}\]
(here $h\otimes\Omega=h$ by convention). Then $\ell(h)^{*}$ exists
and is given by\[
\ell(h)^{*}h_{1}\otimes\cdots\otimes h_{n}=\langle h,h_{1}\rangle h_{2}\otimes\cdots\otimes h_{n}\]
and $\ell(h)^{*}\Omega=0$. The operator $\ell(h)^{*}$ is also called
the annihilation operator.

These operators satisfy\[
\ell^{*}(h)\ell(g)=\langle h,g\rangle1.\]
In particular, the map\[
h\mapsto\ell(h)\]
 is a linear isometry between $H$ (with its Hilbert space norm) and
the closed linear span of $\{\ell(h):h\in H\}$, taken with the operator
norm.

\subsubsection{Relation with non-crossing diagrams.}

Let $h_{1},\ldots,h_{n}\in H$ be an orthonormal family. Let $\ell_{j}=\ell(h_{j})$.
Thus $\ell_{i}^{*}\ell_{j}=\delta_{ij}1$.

The joint distribution of the family $\{\ell_{1},\ell_{1}^{*},\ldots,\ell_{n},\ell_{n}^{*}\}$
(also known as the $*$-distribution of $\{\ell_{1},\ldots,\ell_{n}\}$)
has a nice combinatorial description.

Suppose that we are interested in\[
\phi(\ell_{i(1)}^{g(1)}\cdots\ell_{i(k)}^{g(k)}),\]
where $i(j)\in\{1,\ldots,n\}$ and $g(j)\in\{\cdot,*\}$, $j=1,\ldots,k$
(by $\ell_{j}^{g}$ we mean $\ell_{j}^{*}$ if $g=*$ and $\ell_{j}$
if $g=\cdot$). 

Mark $k$ points on the $x$-axis in half-plane $\{(x,y):y\geq0\}$
at positions $(1,0),\ldots,(k,0)$, and color them by $n$ colors,
so that the $j$-th point point is colored with the $i(j)$-th color.
Attach to the $j$-th point the line segment from $(j,0)$ to $(j,1)$.
Orient this segment upwards (towards infinity) if $g(j)=\cdot$ and
orient it downwards (toward the $x$-axis) if $g(j)=*$. Color the
segment the same way as the $j$-th point, from which it is drawn.

Then there exists at most one way of drawing a diagram so that:

\begin{itemize}
\item The upper end of every segment is connected to the upper end of exactly
one other segment, and all segments connected together have the same
color;
\item Orient each line connecting two segments counter-clockwise. Then the
orientation of the line is compatible with the orientation of the
segments;
\item The lines do not cross.
\end{itemize}
It is not hard to prove that $\phi(\ell_{i(1)}^{g(1)}\cdots\ell_{i(k)}^{g(k)})=1$
iff such a diagram exists, while $\phi(\ell_{i(1)}^{g(1)}\cdots\ell_{i(k)}^{g(k)})=0$
otherwise.

\subsubsection{Moments of $\ell_{1}+\ell_{1}^{*}$.\label{sub:semicirc}}

Utilizing this description one can prove, for example, that\[
\phi((\ell_{1}+\ell_{1}^{*})^{k})=C_{k},\]
where $C_{k}$ is the number of non-crossing pairings between the
integers $\{1,\ldots,k\}$. Recall that a pairing of $\{1,\ldots,k\}$
is an equivalent relation on this set, so that each equivalence class
has exactly two elements. Non-crossing pairings are ones for which
one can draw lines above the real axis $\mathbb{R}\supset\{1,\ldots,k\}$,
connecting the equivalent classes of the pairings, and having no intersections
(more generally, one can in a similar way define non-crossing partitions
of the set $\{1,\ldots,k\}$). We shall later see that the moments
of $X=(\ell_{1}+\ell_{1}^{*})$ are related to the semicircle law.

Non-crossing diagrams and non-crossing partitions have a very deep
connection with free probability; this connection is beyond the scope
of these notes (see e.g. \cite{speicher:thesis}). We will point out
later, however, how this connection explains the relationship between
freeness and large random matrices.

\subsubsection{Free independence.}

Let $A=C^{*}(\ell(h):h\in H)$ and let $\phi:A\to\mathbb{C}$ be given
by\[
\phi(a)=\langle\Omega,a\Omega\rangle.\]
The $C^{*}$-algebra $A$ is an extension of the Cuntz algebra $O_{n}$,
$n=\dim H$ if $n<\infty$ and is isomorphic to $O_{\infty}$ if $\dim H=\infty$.

It is not hard to prove that if $H_{1}\perp H_{2}$ are two subspaces
of $H$, then the algebras\[
C^{*}(\ell(h):h\in H_{1})\ \textrm{and}\  C^{*}(\ell(h):h\in H_{2})\]
are freely independent in $(A,\phi)$.

\subsection{Free Central Limit Theorem.}

\subsubsection{Convergence in moments. }

We say that a sequence $X_{n}$ of random variables converges in moments
to the law a random variable $X$, if $\mu_{X_{n}}\to\mu_{X}$ in
moments; that is to say, for any $p\geq0$,\[
E(X_{n}^{p})\to E(X^{p}).\]
This definition makes sense verbatim (with the replacement of $E$
by $\phi$) in the setting of a non-commutative probability space.

\subsubsection{Classical CLT}

Let $X_{1},\ldots,X_{n},\ldots$ be independent random variables,
so that for all $j$, $E(X_{j})=0$, $E(X_{j}^{2})=1$, and so that
for any $p\geq0$, $\sup_{n}E(X_{n}^{p})\leq C_{p}$ for some constants
$C_{p}<\infty$. The classical central limit theorem states:

\begin{thm}
Let\[
Y_{n}=\frac{1}{\sqrt{n}}(X_{1}+\cdots+X_{n}).\]
Then the laws of the random variables $Y_{n}$ converge in moments
to the Gaussian law $\mu_{\operatorname{Gauss}}$ given by\[
d\mu_{\operatorname{Gauss}}(t)=\frac{1}{\sqrt{2\pi}}\exp(-t^{2}/2)dt.\]

\end{thm}
The main tool used in the proof of this theorem is the fact that if
$Z_{1}$ and $Z_{2}$ are independent random variables, then the law
of their sum is given by a convolution formula:\[
\mu_{Z_{1}+Z_{2}}=\mu_{Z_{1}}*\mu_{Z_{2}}.\]
One then utilizes the fact that the Fourier transform $\hat{\cdot}$
satisfies\[
\widehat{\mu*\nu}=\hat{\mu}\cdot\hat{\nu}.\]
Thus if we write $L_{\mu}=\log\hat{\mu},$ then\[
L_{\mu_{Z_{1}+Z_{2}}}=L_{\mu_{Z_{1}}}+L_{\mu_{Z_{2}}}.\]
Using this one can compute $L_{\mu_{Y_{n}}}$ and argue that it is
quadratic in $t$. This implies that $\mu_{Y_{n}}$ converge in \emph{moments}
to a measure whose Fourier transform is proportional to $\exp(-t^{2}/2)dt$,
so that $\mu_{Y_{n}}\to\mu_{\operatorname{Gauss}}$.

\subsubsection{Free CLT}

Amazingly, the statement of the free central limit theorem is essentially
the same as that of the classical one. The only difference is the
replacement of the requirement of independence by that of free independence.
This is only a single example of a surprising number of parallels
between the behavior of independent and freely independent random
variables.

Let $X_{1},\ldots,X_{n},\ldots$ be freely independent random variables,
so that for all $j$, $\phi(X_{j})=0$, $\phi(X_{j}^{2})=1$, and
so that for any $p\geq0$, $\sup_{n}\phi(|X_{n}^{p}|)\leq C_{p}$
for some constants $C_{p}<\infty$. The classical central limit then
states:

\begin{thm}
\cite{DVV:free} Let\[
Y_{n}=\frac{1}{\sqrt{n}}(X_{1}+\cdots+X_{n}).\]
Then the laws of the random variables $Y_{n}$ converge in moments
to the Gaussian law $\mu_{\operatorname{semicirc}}$ given by\[
d\mu_{\operatorname{semicirc}}(t)=\frac{1}{2\pi}\sqrt{4-t^{2}}dt.\]

\end{thm}
We will postpone the proof of this theorem until we get to talk about
the $R$-transform. For now let us just note that we need a tool to
compute the distribution of $Z_{1}+Z_{2}$ in terms of the distributions
of $Z_{1}$ and $Z_{2}$ if $Z_{1}$ and $Z_{2}$ are freely independent.

\subsection{Free Harmonic Analysis.}

The corresponding classical problem was involved computing the convolution
of two measures via the Fourier transform.

\subsubsection{Free additive convolution.}

By analogy with the classical situation, Voiculescu gave the following
definition:

\begin{defn}
\cite{DVV:free} Let $\mu_{1}$ and $\mu_{2}$ be two probability
measures on $\mathbb{R}$. We define their free additive convolution
$\mu_{1}\boxplus\mu_{2}$ to be the law of the random variable $Z_{1}+Z_{2}$,
where $Z_{1}$ and $Z_{2}$ are freely independent, and $\mu_{Z_{j}}=\mu_{j}$,
$j=1,2$.
\end{defn}
Since $Z_{1},Z_{2}$ are free in $(A,\phi)$, the freeness condition
determines the restriction of $\phi$ to $\textrm{Alg}(Z_{1},Z_{2})$
in terms of the restrictions of $\phi$ to $\textrm{Alg}(Z_{j})$,
$j=1,2$. Thus the joint distribution of $Z_{1}$ and $Z_{2}$ depends
only on $\mu_{Z_{1}}$ and $\mu_{Z_{2}}$. Thus the distribution of
$Z_{1}+Z_{2}$ (which depends only on the joint distribution of $Z_{1}$
and $Z_{2}$) depends only on $\mu_{Z_{1}}=\mu_{1}$ and $\mu_{Z_{2}}=\mu_{2}$.
It follows that $\mu_{1}\boxplus\mu_{2}$ is well-defined.

Note that $\boxplus$ is an operation on the space of probability
measures on $\mathbb{R}$.

\begin{example}
Let $\mu$ be a probability measure and let $\delta_{x}$ be the point
mass at $x$. Then $\mu\boxplus\delta_{x}=\mu_{x}$, the translate
of $\mu$ by $x$. In particular, $\mu\boxplus\delta_{x}$ is the
same as the classical convolution $\mu*\delta_{x}$.
\end{example}

\subsubsection{$R$-transform.}

There is a free analog of the logarithm of the Fourier transform,
which linearizes free additive convolution.

Let $\mu$ be a probability measure on $\mathbb{R}$, and let\[
G_{\mu}(\zeta)=\int_{\mathbb{R}}\frac{d\mu(t)}{\zeta-t},\qquad\Im\zeta>0\]
be a function defined in the upper half-plane. This function is sometimes
callled the Cauchy transform of $\mu$.

If $\mu$ has moments of all orders (e.g., if it is compactly supported),
$G_{\mu}$ is a power series in $1/\zeta$, and we have\[
G_{\mu}(\zeta)=\frac{1}{\zeta}\sum_{p\geq0}\mu_{p}\zeta^{-p},\]
where\[
\mu_{p}=\int_{\mathbb{R}}t^{p}d\mu(t)\]
are the moments of $\mu$. Thus $G_{\mu}$ is the generating function
for the moments of $\mu$.

Define $R_{\mu}(z)$ by the equation\[
G_{\mu}\left(\frac{1}{z}+R_{\mu}(z)\right)=z.\]
It turns out that $R_{\mu}(z)$ is analytic in a certain region in
$\mathbb{C}$; however, one can simply understand it as a formal power
series in $z$ and regard the equation above as an equation involving
composition of formal power series.

Voiculescu proved the following linearization theorem, which shows
that the map $\mu\mapsto R_{\mu}$ is a free analog of the logarithm
of the Fourier transform.

\begin{thm}
\label{thm:RTransf}\cite{DVV:free} Let $R_{\mu}(z)=\sum_{n\geq0}\alpha_{n+1}z^{n}$
be the $R$-transform of $\mu$. Then:\\
(a) $\alpha_{n}$ is a universal polynomial expression in the first
$n$ moments of $\mu$;\\
(b) $R_{\mu}(z)=z$ if and only if $\mu=\mu_{\operatorname{semicirc}}$;
i.e., $d\mu(t)=\frac{1}{2\pi}\sqrt{4-t^{2}}dt$;\\
(c) $R_{\mu_{1}\boxplus\mu_{2}}(z)=R_{\mu_{1}}(z)+R_{\mu_{2}}(z)$;\\
(d) If $Y$ has law $\mu$ and $\lambda\in\mathbb{R}$, then $R_{\mu_{\lambda Y}}(z)=\lambda R_{\mu}(\lambda z)$.
\end{thm}

\subsubsection{Proof of additivity of $R$-transform.}

We will sketch a proof of (a), (b) and (c). We start with a Lemma.

\begin{lem}
\label{lem:combRTransf}Let $X\in(M,\psi)$ be a non-commutative random
variable. Fix $h\in\mathbb{C}$, $\Vert h\Vert=1$. For a sequence
of numbers $a_{1},a_{2},\ldots$, let\[
Y_{N}=Y_{N}^{\{ a_{j}\}_{j=1}^{\infty}}=\ell_{1}^{*}+\sum_{j=0}^{N}a_{j+1}\ell_{1}^{j}\in(C^{*}(\ell(\mathbb{C})),\phi)\]
acting on the full Fock space $F(\mathbb{C})$. Then there exists
a unique sequence of numbers $a_{1},a_{2},\ldots,$, so that for each
$N$,\[
\psi(X^{j})=\phi(Y_{N}^{j}),\quad\forall0\leq j\leq N+1.\]
Moreover, each $a_{k+1}$ is a polynomial in $\{\psi(X^{j}),0\leq j\leq k+1\}$,
and this polynomial is universal, and does not depend on $X$.
\end{lem}
The proof is based on an inductive argument and the combinatorial
formula for moments of free creation operators.

\subsubsection{Combinatorial definition of $R$-transform.}

Given $X$, let $a_{1},a_{2},\ldots$ be as in the Lemma above. Consider
the formal power series\[
R_{\mu}(z)=\sum_{n\geq0}a_{n+1}z^{n}.\]
For now we'll consider $R_{\mu}$ given by this new definition, and
call it the {}``combinatorial $R$-transform''. We shall later prove
that $R_{\mu}(z)$ satisfies our old analytic definition in terms
of $G_{\mu}$ given above; in particular, it will follow that $\alpha_{n}=a_{n}$.

\subsubsection{Additivity of combinatorial $R$-transform.}

\begin{prop}
$R_{\mu_{1}\boxplus\mu_{2}}=R_{\mu_{1}}+R_{\mu_{2}}$.
\end{prop}
\begin{proof}
Let $\ell_{1},\ell_{2}$ be two free creation operators on the free
Fock space $F(\mathbb{C}^{2})$, associated to a pair of orthonormal
vectors.

Given $\mu_{1}$ and $\mu_{2}$, let\[
Y_{1}(n)=\ell_{1}^{*}+\sum_{k\leq n}a_{k+1}\ell_{1}^{k},\quad Y_{2}(n)=\ell_{2}^{*}+\sum_{k\leq n}b_{k+1}\ell_{2}^{k}\]
be random variables in $C^{*}(\ell_{1})$, $C^{*}(\ell_{2})$, respectively,
so that their first $n$ moments are the same as the first $n$ moments
of $\mu_{1}$ and $\mu_{2}$, respectively.

Since $C^{*}(\ell_{1})$ and $C^{*}(\ell_{2})$ are freely independent,
$Y_{1}(n)$ and $Y_{2}(n)$ are freely independent. Since moments
of order up to $n$ of $Y_{1}(n)+Y_{2}(n)$ depend only on the moments
of order up to $n$ of $Y_{1}(n)$ and $Y_{2}(n)$, we see that the
moments of order up to $n$ of $\mu_{1}\boxplus\mu_{2}$ and $Y_{1}(n)+Y_{2}(n)$
are the same.

We leave to the reader the combinatorial exercise to check that the
moments of $Y_{1}(n)+Y_{2}(n)$ are the same as the moments of\[
Y_{3}(n)=\ell_{3}^{*}+\sum(a_{k+1}+b_{k+1})\ell_{3}^{k}.\]
By the uniqueness statement in Lemma \ref{lem:combRTransf}, it follows
that $R_{\mu_{1}\boxplus\mu_{2}}=R_{\mu_{1}}+R_{\mu_{2}}$ as claimed.
\end{proof}

\subsubsection{Analytic and combinatorial $R$-transforms are the same.}

It now remains to prove that the combinatorial $R$-transform $R_{\mu}(z)=\sum a_{n+1}z^{n}$
satisfies the formula relating it to the Cauchy transform $G_{\mu}$
(and so $\alpha_{n}=a_{n}$). The proof of the following proposition
is due to Haagerup \cite{haagerup:r-transform}.

\begin{prop}
Let $K_{\mu}(z)=\frac{1}{z}+R_{\mu}(z)=z^{-1}+\sum a_{k+1}z^{k}$.
With the above notation, one has\[
G_{\mu}(K_{\mu}(z))=z\]
and\[
K_{\mu}(G_{\mu}(\zeta))=\zeta,\]
both equalities interpreted in terms of composition of formal power
series.
\end{prop}
\begin{proof}
Let $\ell$ be a free creation operator corresponding to a unit vector
$e$, and acting on the full Fock space $F(\mathbb{C})$. Let $x=\ell^{*}+f(\ell)$,
where $f$ is a polynomial with real coefficients. Thus by definition,
$R_{\mu_{x}}(z)=f(z)$.

For $z\in\mathbb{C}$ with $|z|<1$, consider the vector\[
\omega_{z}=(1-z\ell)^{-1}\Omega=\Omega+\sum_{n=1}^{\infty}z^{n}e^{\otimes n}.\]
Then\[
\ell\omega_{z}=\sum_{n=0}^{\infty}z^{n}e^{\otimes(n+1)}=\frac{1}{z}(\omega_{z}-\Omega),\quad0<|z|<1.\]
Similarly,\[
\ell^{*}\omega_{z}=\sum_{n=1}^{\infty}z^{n}e^{\otimes(n-1)}=z\omega_{z},\quad|z|<1.\]
Thus $\omega_{z}$ is an eigenvector for $\ell^{*}$ with eigenvalue
$z$. Hence\begin{eqnarray*}
x^{*}\omega_{z} & = & (\ell+f(\ell^{*}))\omega_{z}=\ell\omega_{z}+f(z)\omega_{z}\\
 & = & \frac{1}{z}(\omega_{z}-\Omega)+f(z)\omega_{z}\\
 & = & (\frac{1}{z}+f(z))\omega_{z}-\frac{1}{z}\Omega,\quad0<|z|<1.\end{eqnarray*}
It follows that\[
\frac{1}{z}\Omega=\left(\left(\frac{1}{z}+f(z)\right)1-x^{*}\right)\omega_{z}.\]
Now choose $0<\delta<1$, so that $\left(\left(\frac{1}{z}+f(z)\right)1-x^{*}\right)$
is invertible for $0<|z|<\delta$. This is possible, since $\lim_{z\to0}\left|\frac{1}{z}+f(z)\right|=\infty$
(since $f(z)$ is a polynomial). Hence\[
\left(\left(\frac{1}{z}+f(z)\right)1-x^{*}\right)^{-1}\Omega=z\omega_{z};\]
thus\begin{eqnarray*}
\phi\left(\left(\left(\frac{1}{z}+f(z)\right)1-x^{*}\right)^{-1}\right) & = & \left\langle \left(\left(\frac{1}{z}+f(z)\right)1-x^{*}\right)^{-1}\Omega,\Omega\right\rangle \\
 & = & z\langle\omega_{z},\Omega\rangle=z.\end{eqnarray*}
Since by definition of $G_{\mu}$,\[
G_{\mu}(\lambda)=\phi\left((\lambda1-x)^{-1}\right)=\overline{\phi\left((\bar{\lambda}1-x^{*})^{-1}\right)}.\]
Since all of the coefficients of the power series $G_{\mu}(\lambda)$
are real, we get that\[
\overline{G_{\mu}(\lambda)}=G_{\mu}(\bar{\lambda}),\]
so that\[
G_{\mu}(\bar{\lambda})=\phi\left((\bar{\lambda}1-x^{*})^{-1}\right).\]
We now substitute $\bar{\lambda}=\frac{1}{z}+f(z)$ to get\[
G_{\mu}\left(\left(\frac{1}{z}+f(z)\right)\right)=z.\]
 We also see that $G$ is invertible with respect to composition on
some neighborhood. Applying its inverse to both sides, and remembering
that $f(z)=R_{\mu}(z)$, we get that\[
K_{\mu}(z)=R_{\mu}(z)+\frac{1}{z}=G_{\mu}^{-1}(z)\]
as claimed.

This concludes the proof in the case that $R_{\mu}$ is a polynomial;
the general statement can be deduced from this partial case by taking
limits.
\end{proof}
We have thus proved (a) and (c) of Theorem \ref{thm:RTransf}.

\subsubsection{Semicircular variables.}

Let us prove (b). Assume that $R_{\mu}(z)=z$. Then $K_{\mu}(z)=\frac{1}{z}+z$
and\[
\frac{1}{G_{\mu}(\zeta)}+G_{\mu}(\zeta)=\zeta.\]
Solving this gives\[
G_{\mu}(\zeta)=\frac{\zeta-\sqrt{\zeta^{2}-4}}{2}.\]
One can recover $\mu$ from $G_{\mu}$ by the formula\[
d\mu(t)=\lim_{s\downarrow0}\frac{1}{\pi}G_{\mu}(t+is)dt.\]
Since $G_{\mu}(\zeta)\to0$ as $\zeta\to\infty$ (as is apparent from
the integral formula for the Cauchy transform), the branch of the
square root must be chosen so that $\sqrt{\zeta^{2}-4}>0$ for $\zeta$
real and large. It follows that\[
d\mu(t)=\frac{1}{2\pi}\sqrt{4-t^{2}}dt,\quad t\in[-2,2],\]
and $d\mu(t)=0$ outside of this interval. 

Note that if $R_{\mu}(z)=z$, then $\alpha_{j}=a_{j}=0$ unless $j=2$.
It follows that the variable $\ell_{1}+\ell_{1}^{*}$ on the Fock
space $F(\mathbb{C})$ has semicircular distribution.

\subsubsection{Proof of free CLT}

We are now ready to give a proof of the free central limit theorem.

Let $X_{1},\ldots,X_{n},\ldots$ be freely independent random variables
satisfying the assumptions of the free central limit theorem, and
let\begin{eqnarray*}
Z_{n} & = & (X_{1}+\cdots+X_{n}),\\
Y_{n} & = & \frac{1}{\sqrt{n}}(X_{1}+\cdots+X_{n})=\frac{1}{\sqrt{n}}Z_{n}.\end{eqnarray*}
Let $\nu_{n}$ be the law of $X_{n}$, let $\mu_{n}$ be the law of
$Y_{n}$ and let $\lambda_{n}$ be the law of $Z_{n}$. Thus\[
\lambda_{n}=\nu_{1}\boxplus\cdots\boxplus\nu_{n}\]
and because of additivity of $R$-transform,\[
R_{\lambda_{n}}(z)=R_{\nu_{1}}(z)+\cdots+R_{\nu_{n}}(z).\]
Write $R_{\lambda_{n}}(z)=\sum_{p}\alpha_{p+1}^{(n)}z^{p}$. Since
the coefficient of $z^{p}$ in $R_{\nu_{j}}(z)$ is a universal polynomial
in the moments up to order $p$ of $X_{j}$, and $\sup_{j}|\phi(X_{j}^{p})|<\infty$,
it follows that $|\alpha_{p+1}^{(n)}|\leq n\cdot K_{p}$, where $K_{p}$
are some constants independent of $n$. 

Thus\[
R_{\mu_{n}}(z)=\frac{1}{\sqrt{n}}R_{\lambda_{n}}\left(\frac{z}{\sqrt{n}}\right)=\sum_{p}\frac{\alpha_{p+1}^{(n)}}{n^{\frac{p+1}{2}}}z^{p}.\]
If $p>1$ is fixed, the estimate\[
|\alpha_{p+1}^{(n)}|\leq nK_{p}\]
implies that\[
\frac{\alpha_{p+1}^{(n)}}{n^{\frac{p+1}{2}}}\to0.\]
If $p=0$, the fact that $\phi(X_{n})=0$ so that $\phi(Y_{n})=0$
implies that $\alpha_{1}^{(n)}=0$ for all $n$. Finally, the fact
that $\phi(X_{n}^{2})=1$ implies that $\phi(Y_{n}^{2})=1$ and $\alpha_{2}^{(n)}=1$
for all $n$.

We conclude that\[
R_{\mu_{n}}(z)\to z\]
 as $n\to\infty$ in the sense of coefficient-wise convergence of
formal power series. Since the $p$-th moment of $\mu_{n}$ is a universal
polynomial in the first $p$ coefficients of the power series $R_{\mu_{n}}(z)$,
it follows that the $p$-th moment of $\mu_{n}$ converges to the
$p$-th moment of the unique measure $\mu$ for which $R_{\mu}(z)=z$.
We saw above that this implies that $\mu$ is then the semicircle
measure, and so $\mu_{n}\to\mu_{\operatorname{semicirc}}$.

\subsection{Further topics.}

We already briefly touched upon the amazing correspondence between
various theorems in the classical and free context. There are several
other instances of this. For example, one can consider the free analog
of infinite divisibility. A measure $\mu$ is called infinitely divisible
if for any $n$ there is a measure $\mu_{n}$ so that $\mu$ is the
$n$-fold convolution $\mu_{n}*\cdots*\mu_{n}$. One can say that
$\mu$ is freely infinitely divisible if for each $n$ there is a
measure $\mu_{n}$ so that $\mu$ is the $n$-fold \emph{free} convolution
$\mu_{n}\boxplus\cdots\boxplus\mu_{n}$. Remarkably, there is a one-to-one
correspondence between the classically infinitely divisible measures
and the free ones. A similar situation occurs when considering stable
and freely stable laws. 

There is a also a notion of \emph{multiplicative} free convolution,
based on taking products of non-commutative random variables.

The reader is encouraged to consult \cite{dvv:lectures} for more
details.

\section{Random Matrices and Free Probability.}

One of the most important advances in free probability theory was
Voiculescu's discovery that free probability theory describes the
asymptotic distribution of certain large random matrices. This has
led to a number of applications of free probability theory, both to
spectral computations for random matrices, and to von Neumann algebras.
The latter applications rely on the somewhat unexpected presence of
a {}``matricial'' structure in free probability theory: if one takes
several square arrays of certain free random variables and creates
several matrices out of these arrays, then the resulting matrices
have surprising freeness properties (for example, the resulting matrices
may be freely independent).

\subsection{Random matrices.}

A \emph{random matrix} is a matrix, whose entries are random variables.
One can also think of a random matrix as a matrix-valued random variable,
i.e., as a randomly chosen matrix. Any Borel function of a random
matrix becomes then a random variable. For example, the eigenvalues
of a random matrix (being functions of its entries) are themselves
random variables.

\subsubsection{Expected distributions.}

Let $X_{N}$ be a self-adjoint random matrix of size $N\times N$.
We think of $X_{N}$ as a function $X_{N}:\Sigma\to M_{N}(\mathbb{C})$
on some probability space $(\Sigma,\sigma)$. Integration with respect
to $\sigma$ has the meaning of taking the expected value and will
be denoted by $E$.

One is frequently interested in the expected proportion of the eigenvalues
of $X_{N}$ that lie in a given interval $[a,b]$:\begin{eqnarray*}
\Lambda_{N}([a,b]) & = & \frac{1}{N}\textrm{Expected }\#\{\textrm{eigenvalues of }X_{N}\textrm{ in }[a,b]\}\\
 & = & E(\frac{1}{N}\#\{\textrm{eigenvalues of }X_{N}(t)\ \textrm{in }[a,b]\})\end{eqnarray*}
Let $\lambda_{1}(t),\ldots,\lambda_{N}(t)$ be the eigenvalues of
$X(t)$, listed with multiplicity, and viewed as random variables.
Let\[
\nu_{N}^{t}=\frac{1}{N}\sum_{j=1}^{N}\delta_{\lambda_{j}(t)}\]
 be a random measure associated with this list of eigenvalues (we
say that $\nu_{N}^{t}$ is random to emphasize that it depends on
$t$, i.e., is a measure-valued random variable). Then\[
\Lambda_{N}([a,b])=E(\nu_{N}^{t}([a,b]))\]
is the expected value of $\nu_{N}$. Thus if we set\[
\mu_{N}=E(\nu_{N}^{t})\]
we obtain that\[
\Lambda_{N}([a,b])=\mu_{N}([a,b]).\]
Note that\[
\nu_{N}^{t}=\frac{1}{n}\textrm{Tr}\circ\sigma_{N}^{t},\]
where $\sigma_{N}^{t}$ is the spectral measure of $X_{N}(t)$. In
other words, $\nu_{N}^{t}$ is the distribution of $X_{N}^{t}$, when
viewed as a random variable in $(M_{N}(\mathbb{C}),\frac{1}{N}\textrm{Tr})$.
Thus $\mu_{N}$ is the {}``expected value of the distribution of
$X_{N}$''.

\subsection{Asymptotics of random matrices.}

We are mainly interested in the asymptotics of the expected number
of eigenvalues of a random matrix in a given interval. In other words,
we are interested in studying the asymptotics of the measure $\mu_{N}$
as $N\to\infty$.

It should be mentioned that the eigenvalue distributions of random
matrices have been studied in several ways. Instead of looking at
the expected numbers of eigenvalues, there is also interest in the
behavior of eigenvalue spacing (normalized so that the average spacing
is $1$). One is also interested in the behavior of the largest and
smallest eigenvalues (this translates into considering the expected
value of the spectral radius, or the operator norm, of the matrix
$X_{N}$). We have already heard in this conference of the significant
progress recently made by Haagerup and Thorbjornsen on the latter
problem, in the case that $X_{N}$ is an arbitrary polynomial of a
$k$-tuple of Gaussian random matrices.

\subsubsection{Wigner's theorem for Gaussian random matrices.}

Let $X_{N}$ be a self-adjoint random matrix, whose entries are $g_{ij}$,
$1\leq i,j\leq N$, determined as follows. The variables $\{ g_{ij}:i\leq j\}$
are independent; if $i<j$, then $g_{ij}$ is a centered complex Gaussian
random variable of variance $\frac{1}{N}$; if $i=j$, then $g_{ij}$
is a centered real Gaussian random variable of variance $\frac{2}{N}$.
Finally, if $i>j$, $g_{ij}=\overline{g_{ji}}$.

One can think of the random matrix $X_{N}$ as a map\[
X_{N}:(\Sigma,\sigma)\to M_{N}(\mathbb{C}).\]
Here $\Sigma=M_{N}(\mathbb{C})$ is the space of complex $N\times N$
matrices, $X_{N}$ is the map\[
A\mapsto\frac{A+A^{*}}{2},\]
and $\sigma$ is the Gaussian measure on $\Sigma$ given by\[
d\sigma(A)=\alpha_{N}e^{-\frac{1}{N}\textrm{Tr}(A^{*}A)}dA,\]
for a suitable constant $\alpha_{N}$.

Let $\mu_{N}$ be as before the expected value of the distribution
of $X_{N}$. Then $\mu_{N}\to\mu_{\operatorname{semicirc}}$ weakly
as $N\to\infty$. This is a very old result, going back to the work
of Wigner in 1950s \cite{wigner}.

It turns out that the semicircle law is fairly universal for matrices
with independent identically distributed entries. In fact, Wigner's
original work involved matrices $X_{N}$ whose entries were not Gaussian,
but random signs.

\subsubsection{Voiculescu's asymptotic freeness results.}

The semicircular law also arose in free probability theory as the
central limit law. Voiculescu showed that this is not just a coincidence:
families of certain $N\times N$ random matrices behave as free random
variables in the large $N$ asymptotics.

For each $N$, let $D_{N}$ be a diagonal matrix; assume that the
operator norms $\Vert D_{N}\Vert$ are uniformly bounded in $N$,
and assume that the distribution of $D_{N}$ (as an element of $(M_{N}(\mathbb{C}),\frac{1}{N}\textrm{Tr})$)
converges in moments to a limit measure $\nu$. Let $X_{N}^{(1)},\ldots,X_{N}^{(k)}$
be random matrices described as follows. Let $\Sigma=M_{N}(\mathbb{C})^{k}$
with the measure $\sigma$ given by\[
d\sigma(A_{1},\ldots,A_{k})=C_{N,k}e^{-\frac{1}{N}\textrm{Tr}(A_{1}^{*}A_{1}+\cdots+A_{k}^{*}A_{k})}dA_{1}\cdots dA_{k},\]
for a suitable constant $C_{N,k}$. Then $X_{N}^{(p)}$ is the map\[
X_{N}^{(p)}:(A_{1},\ldots,A_{k})\mapsto\frac{A_{p}+A_{p}^{*}}{2}.\]
More explicitly, if we denote by $g_{ij}^{(p)}$ the $i,j$-th entry
of $X_{N}^{(p)}$, then $\{ g_{ij}^{(p)}:1\leq i\leq j\leq N,\ 1\leq p\leq k\}$
form a family of independent centered Gaussian random variables, so
that: $g_{ij}^{(p)}$ is a complex Gaussian of variance $\frac{1}{N}$
if $i<j$; $g_{ii}^{(p)}$ is real Gaussian of variance $\frac{2}{N}$;
and $g_{ij}^{(p)}=\overline{g_{ji}^{(p)}}$ if $i>j$.

The family $(X_{N}^{(1)},\ldots,X_{N}^{(k)})$ is sometimes called
the Gaussian Unitary Ensemble (or GUE) because of the obvious invariance
of their joint distribution under conjugation by $k$ unitaries.

Let $\mu_{N}$ be the distribution of the family $(D_{N},X_{N}^{(1)},\ldots,X_{N}^{(k)})$,
viewed as a linear functional on the space of polynomials in $k+1$
indeterminates.

Then Voiculescu proved:

\begin{thm}
\label{thm:randomMatrix}\cite{DVV:random} Let $(d,x_{1},\ldots,x_{k})$
be a family of free random variables in a non-commutative probability
space $(A,\phi)$, so that $d$ has distribution $\nu$, and $x_{1},\ldots,x_{k}$
have semicircular distribution. Let $\mu$ be the distribution of
this family, and let $\mu_{N}$ be the distribution of $(D_{N},X_{N}^{(1)},\ldots,X_{N}^{(k)})$
as described above. Then as $N\to\infty$, $\mu_{N}\to\mu$ in moments. 

In other words, for any $t$ and any $j_{1},\ldots,j_{t}\in\{1,\ldots,k\}$,
$n_{0},\ldots,n_{t}\in\{0,1,2,\ldots\}$ one has\begin{eqnarray*}
\lim_{N\to\infty}E\left(\frac{1}{N}\operatorname{Tr}(D_{N}^{n_{0}}X_{N}^{(j_{1})}D_{N}^{n_{1}}\cdots X_{N}^{(j_{t})}D_{N}^{n_{t}})\right)\\
=\phi(d^{n_{0}}x_{j_{1}}d^{n_{1}}\cdots x_{j_{t}}d^{n_{t}}).\end{eqnarray*}

\end{thm}
Note that in particular we have that $D_{N}$ and $X_{N}^{(1)},\ldots,X_{N}^{(k)}$
are asymptotically free. One also recovers Wigner's result, since
in particular $\mu_{X_{N}^{(1)}}\to\mu_{x_{1}}$, and $\mu_{x_{1}}$
is the semicircle law.

\subsubsection{Some remarks on the proof.}

We will not prove this theorem here; see e.g. \cite{DVV:book} for
a proof. We shall only sketch the essential combinatorial trick used
in the proof and explain its connection to non-crossing partitions.

We concentrate on the case of a single random matrix $X_{N}$ with
Gaussian entries $g_{ij}$ (depending on $N$).

Consider the value of the moment\begin{eqnarray}
\frac{1}{N}E(\textrm{Tr}(X_{N}^{k})) & = & \frac{1}{N}\sum_{i_{1},\ldots,i_{k}}E(g_{i_{1}i_{2}}g_{i_{2}i_{3}}\cdots g_{i_{k-1}i_{k}}g_{i_{k}i_{1}}).\label{eq:moment}\end{eqnarray}
If $k$ is odd, it is not hard to see that the value of the moment
is zero, so we'll assume that $k$ is even for the remainder of the
proof. 

Since $g_{ij}$ are Gaussian of variance $\frac{1}{N}$, $E(g_{i_{1}i_{2}}g_{i_{2}i_{3}}\cdots g_{i_{k-1}i_{k}}g_{i_{k}i_{1}})$
is zero unless the variable $g_{i_{p}i_{q}}$ entering in the product
{}``pair up'' with another variable $g_{i_{p'}i_{q'}}$ entering
the product, and $i_{p}=i_{q'}$, $i_{q}=i_{p'}$ (so that $g_{i_{p}i_{q}}=\overline{g_{i_{p'}i_{q'}}}$).
That is to say, a term in the sum (\ref{eq:moment}) is zero unless
for some pairing $\pi$ of the set $\{1,\ldots,k\}$ with itself,
the indices $i_{1},\ldots,i_{k}$ satisfy the equations\begin{equation}
i_{s}=i_{r},\  i_{s+1}=i_{r-1}\qquad\textrm{if }s\sim_{\pi}r,\  s\neq r\label{eq:defbypi}\end{equation}
(where $s+1$ is understood as the remainder mod $n$, and $s\sim_{\pi}r$
iff $s$ and $r$ are in the same equivalence class of $\pi$). 

Suppose now that we fix $\pi$ and ask how large a contribution we
can get from all of the terms that satisfy (\ref{eq:defbypi}) for
this given $\pi$. The equations (\ref{eq:defbypi}) can be visualized
as follows. Let $C_{k}$ be the cyclic graph with $k$ edges, numbered
$1$ through $k$. Place $i_{1},\ldots,i_{k}$ on the vertices of
this graph, so that the $j$-th edge, oriented clockwise, has vertices
$i_{j}$ and $i_{j+1}$, in that order ($j+1$ is again understood
modulo $n$). In other words, we can think of the map $j\mapsto i_{j}$
as a function on the vertices of $C_{k}$.

The pairing $\pi$ defines an equivalence relation on the set of edges
of $C_{k}$: edges $r$ and $s$ are equivalent if $r\sim_{\pi}s$.
Form the quotient graph $C_{k}/\sim_{\pi}$ by gluing equivalent edges
with orientation reversed. Then (\ref{eq:defbypi}) is equivalent
to saying that the function $j\mapsto i_{j}$ descends to a function
on the quotient graph $C_{k}/\sim_{\pi}$. The total number of such
functions is $N^{v}$, where $v$ is the number of vertices of $C_{k}/\sim_{\pi}$. 

Because the variance of $g_{ij}$ is $E(g_{ij}\overline{g_{ij}})=\frac{1}{N}$,
we can deduce that the contribution to the sum (\ref{eq:moment})
of those terms that satisfy equations (\ref{eq:defbypi}) for a given
$\pi$ is at most\[
\frac{1}{N}\cdot\left(\frac{1}{N}\right)^{k/2}\cdot N^{v}.\]
The first factor $1/N$ comes from the normalization of the trace;
the term $(1/N)^{k/2}$ comes from bound on the variance; and the
factor $N^{v}$ comes from our estimation of the number of indices
$i_{1},\ldots,i_{k}$ satisfying (\ref{eq:defbypi}). It follows that
the contribution of all of the terms that satisfy (\ref{eq:defbypi})
for a given $\pi$ is negligible (is of order $1/N$) if $v<1+\frac{k}{2}$.

Recall that $C_{k}$ has exactly $k$ edges and that $k$ is even.
Thus $C_{k}/\sim_{\pi}$ has exactly $k/2$ edges. It follows that
$C_{k}/\sim_{\pi}$ has $1+\frac{k}{2}$ vertices exactly if it is
a tree. With a little bit of care, one can show that (\ref{eq:moment})
is then equal to\[
E(\textrm{Tr}(X_{N}^{k}))=\sum_{\pi\textrm{ s.t. }C_{k}/\sim_{\pi}\textrm{ is a tree}}1+O(\frac{1}{N}).\]
On the other hand, we mentioned in \S\ref{sub:semicirc} that the
$k$-th moment of a semicircular element is given by\[
\tau(s^{k})=\sum_{\sigma\in NC(k)}1,\]
where $NC(k)$ stands for the set of non-crossing pairings of $\{1,\ldots,k\}$.
It is not hard to see that if we interpret a pairing $\sigma$ of
$\{1,\ldots,k\}$ as a pairing of edges of $C_{k}$, it is non-crossing
if and only if $C_{k}/\sim_{\sigma}$ is a tree. This concludes the
proof.

\subsection{An application to random matrix theory.}

Keeping the notations of Theorem \ref{thm:randomMatrix}, let $Y_{N}=D_{N}+X_{N}^{(1)}$.
It is not hard to work out the limit distribution of $Y_{N}$ using
free probability tools. Indeed,\[
\mu_{Y_{N}}\to\mu_{d+x_{1}}.\]
On the other hand, $d$ and $x_{1}$ are freely independent. Thus\[
\mu_{d+x_{1}}=\mu_{d}\boxplus\mu_{x_{1}}=\nu\boxplus\mu_{\operatorname{semicirc}}.\]
The computation of the limit distribution of $Y_{N}$ can then be
carried out using the machinery of $R$-transform.

\subsection{Applications to von Neumann algebras.}

Let us say that a non-commutative non-self-adjoint random variable
$y$ is circular if $\Re y$ and $\Im y$ are freely independent and
are semicircular.

If $X_{N}^{(1)}$ and $X_{N}^{(2)}$ are two GUE random matrices,
then $X_{N}^{(1)}+\sqrt{-1}X_{N}^{(2)}$ converges in $*$-distribution
to a circular variable.

If we start with $2n^{2}$ GUE random matrices $X_{N}^{(i,j,1)},$$X_{N}^{(i,j,2)}$,
$1\leq i,j\leq n$, then we can form a new matrix,\[
Y_{N}=(Y_{N}^{(ij)})_{i,j=1}^{n}\]
of size $nN\times nN$, where $Y_{N}^{(ij)}=X_{N}^{(i,j,1)}+\sqrt{-1}X_{N}^{(i,j,2)}$.
It is not hard to see that $\Re\frac{1}{\sqrt{n}}Y_{N},\Im\frac{1}{\sqrt{n}}Y_{N}$
is a pair of GUE random matrices. We thus obtain that $\frac{1}{\sqrt{n}}Y_{N}$
is circular in the limit $N\to\infty$. From this it is not hard to
prove that if $x_{ij}$, $1\leq i,j\leq n$ are a free circular family,
then the matrix\[
y=\frac{1}{\sqrt{n}}(x_{ij})_{i,j=1}^{n}\]
is again circular. In fact, one can use the asymptotic freeness result
to show that if we let $D$ be the algebra of scalar diagonal $n\times n$
matrices, then $d$ is free from $(y,y^{*})$.

This fact underlies the earliest applications of free probability
theory to von Neumann algebras and free group factors. For example,
one has the following result of Voiculescu \cite{DVV:circular}:

\begin{thm}
Let $n$ be an integer, and let $t\in\mathbb{Q}$ be a rational number,
so that $m=\frac{1}{t^{2}}(n-1)+1$ is an integer. Let $p\in L(\mathbb{F}(n))$
be a projection in the free group factor $L(\mathbb{F}(n))$ associated
to the free group on $n$ generators. Assume that $p$ has trace $t$.
Then\begin{equation}
pL(\mathbb{F}(n))p\cong L(\mathbb{F}(m)).\label{eq:compression}\end{equation}

\end{thm}
This theorem has many far-reaching extensions due to Dykema and Radulescu,
see e.g. \cite{DVV:circular,dykema-subfactors,dykema:extmod,dykema:fdim,dykema:interpolated,radulescu2,radulescu:subfact}.
For example, it turns out that it is possible to define for each $t\in(1,+\infty]$
a von Neumann algebra $L(\mathbb{F}(t))$, called an interpolated
free group factor, in such a way that $L(\mathbb{F}(t))$ is the von
Neumann algebra on the free group with $t$ generators, if $t$ is
an integer. Moreover, the compression formula (\ref{eq:compression})
remains valid for non-rational traces of $p$: the result is an interpolated
free group factor with $\frac{1}{t^{2}}(n-1)+1$ generators; the same
formula is valid also for non-integer $n$. 

For a II$_{1}$ factor $N$, its \emph{fundamental group} was defined
by Murray and von Neumann to be the set\[
F(M)=\{\lambda\in(0,+\infty):M\cong pMp\ \textrm{for }p\in M\textrm{ a projection of trace }\lambda\}.\]
Radulescu proved that $F(L(\mathbb{F}(\infty)))=(0,+\infty)$ (Voiculescu's
result quoted above implied that the positive rational numbers $\mathbb{Q}_{+}\subset F(L(\mathbb{F}(\infty)))$).
In fact, it turns out that there is a dichotomy: either all interpolated
free group factors are the same among each other (and also are isomorphic
to $L(\mathbb{F}(\infty))$), and all have $(0,+\infty)$ as their
fundamental groups; or $L(\mathbb{F}(\infty))\not\cong L(\mathbb{F}(t))$
for finite $t$, and $F(L(\mathbb{F}(t)))=\{1\}$ for finite $t$.
It is not known which of the two alternatives holds.

Further developments of these techniques gave information on fundamental
groups of more general free products of von Neumann algebras and on
subfactors of $L(\mathbb{F}(\infty))$ (see e.g. \cite{radulescu:subfact,dykema-subfactors,shlyakht:amalg,shlyakht:semicirc,shlyakht-popa:universal,shlyakht-ueda:subfactors,dykema-radulescu:compressions}).

\section{Free Entropy via Microstates.}

Free entropy was introduced and developed by Voiculescu in a series
of papers \cite{dvv:entropy1,dvv:entropy2,dvv:entropy3,dvv:entropy4,dvv:improvedrandom,dvv:entropy5,dvv:entropy6,dvv:entropySLnZ}
as a free probability analogue of the classical information-theoretic
entropy; see also Voiculescu's survey \cite{dvv:entropysurvey}.

\subsection{Definition of free entropy.}

Voiculescu's original {}``microstates'' approach to free entropy
followed Boltzman's definition of entropy of a macroscopic state.

\subsubsection{Microstates and Macrostates.}

Assume that the macroscopic behavior of a physical system (e.g. gas)
is described by several macroscopic parameters (e.g., pressure, volume
and temperature). Then a \emph{macrostate $s$} is a state of the
system corresponding to certain prescribed values of these parameters. 

Microscopically, the system is made out of a large number of smaller
systems (e.g., the molecules that make up the gas). On this microscopic
level, the system can be described by a \emph{microstate} $s$ that
specifies exactly the states of all of the sub-systems (e.g, the exact
locations and moments of all of the molecules of the gas). If we fix
a macrostate $S$, there are many microstates $s$ that lead to the
same macroscopic state. 

Boltzman's formula is then that the entropy of $S$ must be given
by\[
K\log\#\{ s:\textrm{microscopic state }s\textrm{ leads to macroscopic state }S\}\]
for some constant $K$.

\subsubsection{Matricial microstates.}

Voiculescu's idea is to interpret $x_{1},\ldots,x_{n}\in(A,\phi)$
as a description of a macroscopic state of a system, and as microstates
to take the set of all matrices $X_{1},\ldots,X_{n}$ of a specific
dimension that approximate $x_{1},\ldots,x_{n}$. More precisely for
$x_{1},\ldots,x_{n}$ in a non-commutative probability space $(A,\phi)$,
$x_{j}=x_{j}^{*}$, let $M_{k\times k}^{sa}$ be the space of $k\times k$
self-adjoint matrices, and consider the set\begin{eqnarray*}
 &  & \Gamma(x_{1},\ldots,x_{n};k,l,\varepsilon)=\Big\{(X_{1},\ldots,X_{n})\in(M_{k\times k}^{sa})^{n}:\\
 &  & \textrm{ }\qquad\textrm{for any word }w\ \textrm{in }n\ \textrm{letters of length at most }l,\\
 &  & \qquad\qquad|\frac{1}{N}\textrm{Tr}(w(X_{1},\ldots,X_{n}))-\phi(w(x_{1},\ldots,x_{n}))|<\varepsilon\Big\}.\end{eqnarray*}
In other words, we are considering a weak neighborhood $U$ of the
joint law $\mu_{x_{1},\ldots,x_{n}}$defined by the property that
$\mu'\in U$ iff the value of the law $\mu'$ on all words of length
at most $l$ deviates by no more than $\varepsilon$ from that of
$\mu_{x_{1},\ldots,x_{n}}$. Next, we consider all self-adjoint $k\times k$
matrices $(X_{1},\ldots,X_{n})$ so that\[
\mu_{X_{1},\ldots,X_{n}}\in U.\]
The set $\Gamma(x_{1},\ldots,x_{n};k,l,\varepsilon)$ is called the
set of (matricial) microstates for $x_{1},\ldots,x_{n}$.

\subsubsection{Definition of free entropy.}

Voiculescu then defined the free entropy by\[
\chi(x_{1},\ldots,x_{n})=\inf_{\varepsilon,l}\limsup_{k\to\infty}\frac{1}{k^{2}}\log\textrm{Vol}\Gamma(x_{1},\ldots,x_{n};k,l,\varepsilon)+\frac{n}{2}\log k,\]
where $\textrm{Vol}$ refers to the Euclidean volume associated to
the standard identification of $M_{k\times k}^{sa}$ with $\mathbb{R}^{k^{2}}$.
We use the convention that $\log0=-\infty$.

We should note that $\chi$ depends only on the \emph{law} of $x_{1},\ldots,x_{n}$
and not on the particular realization of this law. It would be also
appropriate to write $\chi(\mu_{x_{1},\ldots,x_{n}})$.

\subsubsection{Relation to Connes' problem.}

Note that there is no a priori reason for $\Gamma(x_{1},\ldots,x_{n};k,l,\varepsilon)$
to be non-empty. Connes has posed a question in \cite{connes:injective}
of whether every II$_{1}$ factor can be embedded into an ultrapower
of the hyperfinite II$_{1}$ factor. It is not hard to see that his
question is equivalent to the question of whether, given $x_{1},\ldots,x_{n}$
in a von Neumann algebra $(A,\phi)$ with $\phi$ a trace, one has
that for any $\varepsilon>0$ and $l>0$ there is a $k$ so that $\Gamma(x_{1},\ldots,x_{n};k,l,\varepsilon)\neq\emptyset$.
This question is open even for $x_{1},\ldots,x_{n}$ elements of the
group algebra of an arbitrary discrete group $\Gamma$.

\subsection{Properties of free entropy.}

Voiculescu gave an explicit formula for the free entropy of a single
variable $x$ with law $\mu$:\[
\chi(x)=\iint\log|s-t|d\mu(s)d\mu(t)+C\]
for a certain universal constant $C$.

Free entropy has a number of nice properties, related to freeness
and analogous to the properties of classical entropy; we list a few,
due to Voiculescu \cite{dvv:entropy2}:

\begin{itemize}
\item If $x_{1},\ldots,x_{n}$ are free, then $\chi(x_{1},\ldots,x_{n})=\chi(x_{1})+\cdots+\chi(x_{n})$.
Furthermore, if $\chi(x_{1},\ldots,x_{n})=\chi(x_{1})+\cdots+\chi(x_{n})\neq-\infty$,
then $x_{1},\ldots,x_{n}$ are freely independent.
\item $\chi(x_{1},\ldots,x_{n},y_{1},\ldots,y_{m})\leq\chi(x_{1},\ldots,x_{n})+\chi(y_{1},\ldots,y_{m})$.
\item $\chi(x_{1},\ldots,x_{n})$ is maximal subject to $\sum\phi(x_{i}^{2})=n^{2}$
iff $x_{1},\ldots,x_{n}$ is a free semicircular family and each $x_{i}$
satisfies $\phi(x_{i}^{2})=1$.
\item If $s_{1},\ldots,s_{n}$ are free semicircular variables, freely independent
from the family $x_{1},\ldots,x_{n}$, then $W^{*}(x_{1},\ldots,x_{n})$
embeds into the ultrapower of the hyperfinite II$_{1}$ factor if
and only if $\chi(x_{1}+\sqrt{\delta}s_{1},\ldots,x_{n}+\sqrt{\delta}x_{n})>-\infty$
for every $\delta>0$. Thus semicircular perturbations (i.e., {}``free
Brownian motion'') have a regularization effect on free entropy.
\end{itemize}
To give but one example of the technical difficulties that working
with $\chi$ presents, one would be able to prove that\[
\chi(x_{1},\ldots,x_{n},y_{1},\ldots,y_{m})=\chi(x_{1},\ldots,x_{n})+\chi(y_{1},\ldots,y_{m})\]
if $(x_{1},\ldots,x_{n})$ and $(y_{1},\ldots,y_{m})$ are free families,
provided that one could argue that the $\limsup$ in the definition
of free entropy is a limit.

\subsubsection{Infinitesimal change of variables formula.}

We end the review of free entropy by mentioning the change of variables
formula \cite{dvv:entropy2}.

Assume that $y_{1},\ldots,y_{n}$ are given as non-commutative power
series in $x_{1},\ldots,x_{n}$: $y_{j}=F_{j}(x_{1},\ldots,x_{n})$.
Assume moreover that the multi-radius of convergence of $F_{j}$ is
large enough to exceed the norms of all $y_{1},\ldots,y_{n}$. Assume
further that $x_{j}=G_{j}(y_{1},\ldots,y_{n})$ for some non-commutative
power series $G_{j}$, and that similarly the multi-radius of convergence
of $G_{j}$ is large enough to exceed the norms of $x_{1},\ldots,x_{n}$. 

Let $M=W^{*}(x_{1},\ldots,x_{n})=W^{*}(y_{1},\ldots,y_{n})$, and
let $\phi$ be the given trace on $M$. Consider the derivation $\partial_{j}:\mathbb{C}[x_{1},\ldots,x_{n}]\to M\bar{\otimes}M$
determined by\[
\partial_{j}(x_{i})=\delta_{ji}1\otimes1.\]
For example,\[
\partial_{2}(x_{1}x_{2}^{2}x_{3}x_{2})=x_{1}\otimes x_{2}x_{3}x_{2}+x_{1}x_{2}\otimes x_{3}x_{2}+x_{1}x_{2}^{2}x_{3}\otimes1.\]
Let\[
J(x_{1},\ldots,x_{n})=\left(J_{ij}(x_{1},\ldots,x_{n})\right)_{ij=1}^{n}\in M_{n}(M\bar{\otimes}M)\]
be the {}``Jacobian'' of $F$: $J_{ij}(x_{1},\ldots,x_{n})=\partial_{i}F_{j}(x_{1},\ldots,x_{n})$.
Then\[
\chi(y_{1},\ldots,y_{n})=\chi(x_{1},\ldots,x_{n})+n\log(|\det|(J(x_{1},\ldots,x_{n}))),\]
where $|\det|$ refers to the Kadison-Fuglede determinant\[
|\det|(J)=\exp(\tau_{M_{n\times n}(M\otimes M)}(\log|J|)).\]
Here $\tau_{M_{n\times n}(M\otimes M)}$ is the tensor product $\frac{1}{n}\textrm{Tr}\otimes\phi\otimes\phi$
of the traces on $M_{n\times n}$ and $M\otimes M$.

The explanation of this formula and the appearance of $J$ is that
the Jacobian of the transformation\[
(X_{1},\ldots,X_{n})\mapsto(F_{1}(X_{1},\ldots,X_{n}),\ldots,F_{n}(X_{1},\ldots,X_{n})),\]
viewed as a map from $M_{k\times k}^{n}\to M_{k\times k}^{n}$ is
naturally a matrix in $M_{n\times n}(\textrm{End}(M_{k\times k}))\cong M_{n\times n}(M_{k\times k}\otimes M_{k\times k})$,
and is given by $J(X_{1},\ldots,X_{n})$.

\subsection{Free entropy dimension.}

Voiculescu's original idea for defining free entropy dimension was
to consider a kind of asymptotic Minkowski dimension of the set of
microstates. We present below an equivalent definition of K. Jung,
which is based on packing dimension instead.

\subsubsection{Packing and covering numbers and Minkowski dimension.}

For a metric space $X$, let $P_{\varepsilon}(X)$ be the \emph{packing
number}~of $X$; that is, the maximal number of disjoint $\varepsilon$-balls
that can be placed inside $X$. Similarly, let $K_{\varepsilon}(X)$
be the \emph{covering} \emph{number}~of $X$; that is, the minimal
number of $\varepsilon$-balls needed to cover $X$.

For a metric space $X$, the \emph{upper uniform packing dimension}
and the \emph{upper uniform covering dimension} are the same and are
defined as\[
\limsup_{\varepsilon\to0}\frac{\log P_{\varepsilon}(X)}{|\log\varepsilon|}=\limsup_{\varepsilon\to0}\frac{\log K_{\varepsilon}(X)}{|\log\varepsilon|}.\]
It is a theorem that if $X\subset\mathbb{R}^{d}$, then both of these
numbers are the same as the \emph{Minkowski dimension} of $X$, which
is given by\[
d-\liminf_{\varepsilon\to0}\frac{\log\textrm{Vol}N_{\varepsilon}(X)}{\log\varepsilon},\]
where $N_{\varepsilon}(X)$ denotes the tubular neighborhood of $X$
of radius $\varepsilon$.

\subsubsection{Free entropy dimension.}

Let $x_{1},\ldots,x_{n}\in(A,\phi)$ be self-adjoint. Then let\begin{eqnarray*}
P_{\delta}(x_{1},\ldots,x_{n}) & = & \inf_{\varepsilon,l}\limsup_{k\to\infty}\frac{1}{k^{2}}\log P_{\delta}(\Gamma(x_{1},\ldots,x_{n};l,k,\varepsilon))\\
K_{\delta}(x_{1},\ldots,x_{n}) & = & \inf_{\varepsilon,l}\limsup_{k\to\infty}\frac{1}{k^{2}}\log K_{\delta}(\Gamma(x_{1},\ldots,x_{n};l,k,\varepsilon)).\end{eqnarray*}
Then K. Jung proved the following theorem \cite{jung:packing}:

\begin{thm}
One has\[
\limsup_{\delta\to0}\frac{P_{\delta}(x_{1},\ldots,x_{n})}{|\log\delta|}=\limsup_{\delta\to0}\frac{K_{\delta}(x_{1},\ldots,x_{n})}{|\log\delta|}.\]
Moreover, if $s_{1},\ldots,s_{n}$ are free semicircular variables,
free from $x_{1},\ldots,x_{n}$, then\[
\limsup_{\delta\to0}\frac{P_{\delta}(x_{1},\ldots,x_{n})}{|\log\delta|}=n-\liminf_{\delta\to0}\frac{\chi(x_{1}^{\delta},\ldots,x_{n}^{\delta}:s_{1},\ldots,s_{n})}{\log\delta^{1/2}},\]
where $x_{j}^{\delta}=x_{j}+\sqrt{\delta}s_{j}$.
\end{thm}
The value of any of these limits is then by definition called the
\emph{free entropy dimension} $\delta_{0}(x_{1},\ldots,x_{n})$.

Here $\chi(x_{1}+\sqrt{\delta}s_{1},\ldots,x_{n}+\sqrt{\delta}s_{n}:s_{1},\ldots,s_{n})$
is the free entropy of $x_{1}+\sqrt{\delta}s_{1},\ldots,x_{n}+\sqrt{\delta}s_{n}$
in the presence of $s_{1},\ldots,s_{n}$; it is a technical modification
of the free entropy $\chi(x_{1}+\sqrt{\delta}s_{1},\ldots,x_{n}+\sqrt{\delta}s_{n})$.
Very roughly, the value of $\chi(x_{1}+\sqrt{\delta}s_{1},\ldots,x_{n}+\sqrt{\delta}s_{n})$
is the asymptotic logarithmic volume of a $\delta^{1/2}$-tubular
neighborhood of the set of microstates for $x_{1},\ldots,x_{n}$.
Thus the number\[
n-\liminf_{\delta\to0}\frac{\chi(x_{1}+\sqrt{\delta}s_{1},\ldots,x_{n}+\sqrt{\delta}s_{n}:s_{1},\ldots,s_{n})}{\log\delta^{1/2}}\]
is a kind of asymptotic Minkowski dimension of the set of microstates.
This was the original definition of free entropy dimension given by
Voiculescu.

We finish this section with an example. 

Let $x_{1},\ldots,x_{n}$ be free semicircular variables. Then $x_{1}+\sqrt{\delta}s_{1},\ldots,x_{n}+\sqrt{\delta}s_{n}$
are also semicircular. In fact\[
\chi(x_{1}+\sqrt{\delta}s_{1},\ldots,x_{n}+\sqrt{\delta}s_{n}:s_{1},\ldots,s_{n})\geq\chi(x_{1},\ldots,x_{n})>-\infty.\]
It follows that $\delta_{0}(x_{1},\ldots,x_{n})=n$. In particular,
the free group factor $L(\mathbb{F}(n))$ can be generated by a family
with free entropy dimension $n$.

\subsection{Properties of free entropy dimension.}

The theory of free entropy dimension has found a number of spectacular
applications to von Neumann algebra theory. For example, Voiculescu
used free entropy dimension to prove that free group factors do not
have Cartan subalgebras; soon thereafter, L. Ge gave a proof that
free group factors are prime, i.e., cannot be written as tensor products
of infinite-dimensional von Neumann algebras. 

One of the main remaining questions about free entropy dimension is
the extent to which $\delta_{0}(x_{1},\ldots,x_{n})$ depends on the
elements $x_{1},\ldots,x_{n}$. Voiculescu asked if $\delta_{0}(x_{1},\ldots,x_{n})$
is an invariant of the von Neumann algebra generated by $x_{1},\ldots,x_{n}$,
taken with a fixed trace. Since $L(\mathbb{F}(n))$ has a generating
family with free entropy dimension equal to $n$, a positive answer
to this question would imply non-isomorphism of free group factors.

\subsubsection{Invariance of $\delta_{0}$.}

Voiculescu proved that $\delta_{0}(x_{1},\ldots,x_{n})$ depends only
on the restriction of the trace to the algebra generated by $x_{1},\ldots,x_{n}$.
In particular, if $\Gamma$ is a discrete group and $x_{1},\ldots,x_{n}\in\mathbb{C}\Gamma$
are self-adjoint generators of the group algebra, then $\delta_{0}(x_{1},\ldots,x_{n})$
depends only on the group. This invariant seems to be related to the
$L^{2}$-cohomology of $\Gamma$; see below.

\subsubsection{Free entropy dimension for a single variable.}

Voiculescu proved that if $X$ has law $\mu$, then \[
\delta_{0}(X)=1-\sum_{t\textrm{ an atom of }\mu}\mu(\{ t\})^{2}.\]
In particular, notice that $\delta_{0}$ is an invariant of the von
Neumann algebra (with a fixed trace) generated by $X$.

\subsubsection{Upper bounds on $\delta_{0}$.}

If $M$ satisfies any of the following conditions, then $\delta_{0}(x_{1},\ldots,x_{n})=1$
for any $x_{1},\ldots,x_{n}\in M$ generating $M$: 

\begin{enumerate}
\item \cite{dvv:entropy3} $M$ has a Cartan subalgebra, i.e., a maximal
abelian subalgebra $A$ so that $M=W^{*}(\{ u\in M\textrm{ unitary}:uAu^{*}=A\})$.
Thus free group factors have no Cartan subalgebras.
\item \cite{dvv:entropy3} $M$ has a diffuse regular hyperfinite subalgebra:
a hyperfinite subalgebra $R$ so that $M=W^{*}(\{ u\in M\textrm{ unitary}:uRu^{*}=R\})$.
This is the case, in particular, if $M=L(\Gamma)$ and $\Gamma$ has
an infinite normal amenable subgroup. Thus free group factors do not
have diffuse regular hyperfinite subalgebras.
\item \cite{dvv:entropy3} $M$ has property $\Gamma$: there is a sequence
of unitaries $u_{n}\in M$, so that $\tau(u_{n})\to0$ but $\Vert u_{n}x-xu_{n}\Vert_{2}\to0$
for all $x\in M$. Free group factors are non-$\Gamma$ by a classical
result of Murray and von Neumann.
\item \cite{ge:entropy2} $M\cong M_{1}\otimes M_{2}$ with $M_{1}$ and
$M_{2}$ infinite-dimensional. Thus free group factors are prime.
\end{enumerate}
In particular, note that $M\not\cong L(\mathbb{F}_{n})*N$ for any
$N$ which be embedded into the ultrapower of the hyperfinite II$_{1}$
factor (e.g., $N=\mathbb{C}$ is already interesting).

There are other conditions assuring upper bounds on $\delta_{0}$;
we mention the work of K. Dykema \cite{dykema:freeenetropy}, M. Stefan
\cite{stephan:thinness} and of Ge and Shen \cite{ge-shen:freeDimPropT}.
Upper estimates on $\delta_{0}$ turned out to be of relevance also
to the theory of type III factors \cite{shlyakht:prime,shlyakht:fullfactor}.

\subsubsection{Lower bounds on $\delta_{0}$.}

K. Jung has proved the following {}``hyperfinite monotonicity result''
\cite{jung-freexentropy}: let $M$ be a diffuse von Neumann algebra,
and assume that $M$ is embeddable in the ultrapower of the hyperfinite
II$_{1}$ factor. Then $\delta_{0}(x_{1},\ldots,x_{n})\geq1$ for
any generators $x_{1},\ldots,x_{n}$.

Combined with the upper estimates, this shows that if $M$ satisfies
any of the properties (1)--(4) above and is embeddable into the ultrapower
of the hyperfinite II$_{1}$ factor, then the value of $\delta_{0}$
is $1$ on any set of generators. In particular, $\delta_{0}$ \emph{is}
an invariant of the entire von Neumann algebra!

Jung has also computed $\delta_{0}$ for arbitrary generators of a
hyperfinite algebra \cite{jung-freexentropy} (which is in general
a direct sum of matrix algebras and a diffuse hyperfinite von Neumann
algebra) and once again found that $\delta_{0}$ is an invariant of
the von Neumann algebra in that case.

\subsection{Relation with $L^{2}$-Betti numbers.}

By \cite{connes-shlyakht:l2betti}, for any generators $(x_{1},\ldots,x_{n})$
of a tracial algebra $(A,\tau)$ one has the inequality relating $\delta_{0}$
to the $L^{2}$-Betti numbers of $A$:\[
\delta_{0}(A)=\delta_{0}(x_{1},\ldots,x_{n})\leq\beta_{1}^{(2)}(A,\tau)-\beta_{0}^{(2)}(A,\tau)+1.\]
In particular, specializing to the case of the group algebra of a
discrete group $\Gamma$, we have that\[
\delta_{0}(\Gamma)\leq b_{1}^{(2)}(\Gamma)-b_{0}^{(2)}(\Gamma)+1,\]
where $b_{j}^{(2)}$ are the $L^{2}$-Betti numbers of the group. 

The same combination of Betti numbers also occurs in Gaboriau's work
on cost of equivalence relations \cite{gaboriau:cost,gaboriau:ell2};
indeed he proves that\[
b_{1}^{(2)}(\Gamma)-b_{0}^{(2)}(\Gamma)+1\leq C(\Gamma),\]
where $C(\Gamma)$ is the cost of $\Gamma$. There are no known examples
in which equality does not hold. 

It is curious that $C(\Gamma)$ measures the {}``optimal number of
generators'' for an equivalence relation induced by $\Gamma$; on
the other hand, $\delta_{0}(x_{1},\ldots,x_{n})$ is known to be $\leq1$
in many cases in which the von Neumann algebra is {}``singly generated''
\cite{popa-ge:thin}.

One obstruction for the equality between $\delta_{0}(\Gamma)$ and
$b_{1}^{(2)}(\Gamma)-b_{0}^{(2)}(\Gamma)+1$ is the fact that the
latter quantity is insensitive to the outcome of Connes' embedding
question (if there is an non-embeddable group, one can manufacture
a non-embeddable group with large Betti numbers by taking free products).

It is also possible to define a {}``relative'' version of Voiculescu's
free entropy dimension for equivalence relations; one can obtain an
invariant of an equivalence relation in this way (see \cite{shlyakht:cost,shlyakht:cost:micro}.

\section{Non-microstates Approach to Free Entropy.}

We have reviewed the microstates definition of free entropy in the
previous lecture. There are several difficulties connected with that
definition. The first is that the involvement of sets of microstates
makes the definition hard to work with technically; as we saw there
are several properties of free entropy (such as additivity for free
families) that one expects to hold, but which one is unable to prove
because of such technical difficulties. Another example of such acute
difficulties arises when one deals with free Fisher information. By
analogy with the classical case, one wants to define the free Fisher
information $\Phi(x_{1},\ldots,x_{n})$ by the formula\[
\Phi(x_{1},\ldots,x_{n})=2\frac{d}{d\varepsilon}\chi(x_{1}+\sqrt{\varepsilon}s_{1},\ldots,x_{n}+\sqrt{\varepsilon}s_{n})\Big|_{\varepsilon=0},\]
where $s_{1},\ldots,s_{n}$ are free semicircular variables, free
from $(x_{1},\ldots,x_{n})$. The definition works fine in the case
that $n=1$ (the explicit formula for $\chi$ is essential), but it
is not clear how to prove that the derivative exists and that the
definition makes sense in the case $n>1$.

The other point is that the definition of the microstates free entropy
subsumes existence of microstates, i.e., embedability into the ultrapower
of the hyperfinite II$_{1}$ factor. A priori, it is not clear why
one should assume this for elements of an arbitrary non-commutative
tracial probability space (although of course if Connes' embedability
question always has an affirmative answer, this second point disappears). 

Voiculescu \cite{dvv:entropy5} gave a new definition of free entropy,
based on an {}``infinitesimal'' approach involving free Fisher information.
This new approach does not involve microstates and for this reason
the resulting entropy bears the name {}``non-microstates'' or {}``microstates-free''.
It is not at present known if the two definitions (microstates and
non-microstates) are the same, except in the one-variable case; and
indeed, showing this would give a positive answer to Connes' embedability
question. Nonetheless, a recent work by Biane, Capitaine and Guionnet
\cite{guionnet-biane-capitaine:largedeviations} shows that the microstates
free entropy is always smaller than the non-microstates entropy. 

To distinguish the two definitions, quantities related to the non-microstates
entropy are denoted by the same letter as their microstates analogs,
but with an asterisk; for example, the non-microstates free entropy
is $\chi^{*}$, and the corresponding free entropy dimension is $\delta^{*}$.

\subsection{A non-rigorous derivation of the non-microstates definition.}

We begin with a (rigorous) consequence of the change of variables
formula for microstates entropy. We shall assume that $x_{1},\ldots,x_{n}$
are in a non-commutative probability space $A$ with a tracial positive
linear functional $\tau$.

\subsubsection{Infinitesimal change of variables.}

Let $P_{1},\ldots,P_{n}$ be polynomials in $n$ indeterminates. Consider
the change of variables\[
x_{j}^{\varepsilon}=x_{j}+\varepsilon P_{j}(x_{1},\ldots,x_{n}).\]
Then for $\varepsilon$ sufficiently small, this change of variables
can be inverted and $x_{j}$ can be expressed as a non-commutative
power series in terms of $y_{1}^{\varepsilon},\ldots,y_{n}^{\varepsilon}$,
so that the multi-radius of convergence of that power series exceeds
the operator norms of $y_{1}^{\varepsilon},\ldots,y_{n}^{\varepsilon}$.
Thus one can apply the change of variables formula and express $\chi(y_{1}^{\varepsilon},\ldots,y_{n}^{\varepsilon})$
in terms of the free entropy $\chi(x_{1},\ldots,x_{n})$ and the logarithm
of the Jacobian of our transformation. Expanding the value of the
logarithm of the Jacobian as a power series in $\varepsilon$ gives
us the \emph{infinitesimal change of variables formula \cite{dvv:entropy4}}:\[
\chi(y_{1}^{\varepsilon},\ldots,y_{n}^{\varepsilon})=\chi(x_{1},\ldots,x_{n})+\varepsilon\sum_{j=1}^{n}\tau\otimes\tau(\partial_{j}P_{j})+O(\varepsilon^{2}).\]

\subsubsection{Conjugate variables.}

Let us now assume that $\partial_{j}:L^{2}(M)\to L^{2}(M\bar{)\otimes}L^{2}(M)$,
with $M=W^{*}(x_{1},\ldots,x_{n})$ has the property that $1\otimes1$
is in the domain of $\partial_{j}^{*}$. Let $\xi_{j}=\partial_{j}^{*}(1\otimes1)\in L^{2}(M)$.
The elements $\xi_{1},\ldots,\xi_{n}$ are called \emph{conjugate
variables} to $(x_{1},\ldots,x_{n})$ and satisfy\[
\langle\xi_{j},Q\rangle=\langle\partial_{j}(Q),1\otimes1\rangle=\tau\otimes\tau(\partial_{j}(Q)),\]
for any polynomial $Q\in\mathbb{C}[X_{1},\ldots,X_{n}]$.

Then our infinitesimal change of variables formula becomes:\[
\chi(y_{1}^{\varepsilon},\ldots,y_{n}^{\varepsilon})=\chi(x_{1},\ldots,x_{n})+\varepsilon\sum_{j=1}^{n}\langle P_{j},\xi_{j}\rangle+O(\varepsilon^{2}).\]

It turns out that conjugate variables are intimately connected with
free Brownian motion. If we let\[
x_{j}^{\varepsilon}=x_{j}+\sqrt{\varepsilon}s_{j},\]
where $s_{j}$ are a free semicircular family, free from $x_{1},\ldots,x_{n}$,
then for any polynomial $Q$ in $n$ indeterminates one can prove
that\[
\tau(Q(x_{1}^{\varepsilon},\ldots,x_{n}^{\varepsilon}))=\tau(Q(x_{1}+\frac{\varepsilon}{2}\xi_{1},\ldots,x_{n}+\frac{\varepsilon}{2}\xi_{n}))+O(\varepsilon^{2}).\]
Thus perturbations by conjugate variables give an {}``approximation''
in law to free Brownian motion; note, however, that while $x_{j}^{\varepsilon}$
no longer lies in $W^{*}(x_{1},\ldots,x_{n})$, $x_{j}+\frac{\varepsilon}{2}\xi_{j}$
does lie in $L^{2}(W^{*}(x_{1},\ldots,x_{n}))$.

Conjugate variables frequently exist. For example, if $x_{1},\ldots,x_{n}$
are a free semicircular family, then $\xi_{j}$ exist and in fact
$\xi_{j}=s_{j}$, $j=1,\ldots,n$. One can show that for any $x_{1},\ldots,x_{n}$
and any $\varepsilon>0$, conjugate variables to the family $(x_{1}+\sqrt{\varepsilon}s_{1},\ldots,x_{n}+\sqrt{\varepsilon}s_{n})$
always exist. In fact, in this case\[
\xi_{j}=E_{W^{*}(x_{1}+\sqrt{\varepsilon}s_{1},\ldots,x_{n}+\sqrt{\varepsilon}s_{n})}\left(\frac{1}{\sqrt{\varepsilon}}s_{j}\right).\]

\subsubsection{Non-rigorous derivation of the formula for $\Phi(x_{1},\ldots,x_{n})$.}

Assume now that $(\xi_{1},\ldots,\xi_{n})$ are conjugate variables
to $(x_{1},\ldots,x_{n})$. Let $(s_{1},\ldots,s_{n})$ be as before
a free semicircular system, free from $(x_{1},\ldots,x_{n})$.

Recall that we want to define the free Fisher information by\[
\Phi(x_{1},\ldots,x_{n})=2\frac{d}{d\varepsilon}\chi(x_{1}+\sqrt{\varepsilon}s_{1},\ldots,x_{n}+\sqrt{\varepsilon}x_{n}).\]
Since $\chi(x_{1},\ldots,x_{n})$ depends only on the law of $x_{1},\ldots,x_{n}$,
and since the laws of $(x_{1}+\sqrt{\varepsilon}s_{1},\ldots,x_{n}+\sqrt{\varepsilon}x_{n})$
and $(x_{1}+\frac{\varepsilon}{2}\xi_{1},\ldots,x_{n}+\frac{\varepsilon}{2}\xi_{n})$
are the same up to higher orders in $\varepsilon$, one would expect
that\[
\Phi(x_{1},\ldots,x_{n})=2\frac{d}{d\varepsilon}\chi(x_{1}+\frac{\varepsilon}{2}\xi_{1},\ldots,x_{n}+\frac{\varepsilon}{2}\xi_{n}).\]
We now assume that $\xi_{j}$ are sufficiently nice functions of $x_{1},\ldots,x_{n}$
so that the infinitesimal change of variables applies. Thus\begin{eqnarray*}
\chi(x_{1}+\frac{\varepsilon}{2},\ldots,x_{n}+\frac{\varepsilon}{2}) & = & \chi(x_{1},\ldots,x_{n})+\varepsilon\sum_{j=1}^{n}\langle\frac{1}{2}\xi_{j},\xi_{j}\rangle+O(\varepsilon^{2})\\
 & = & \chi(x_{1},\ldots,x_{n})+\frac{\varepsilon}{2}\sum_{j=1}^{n}\Vert\xi_{j}\Vert_{L^{2}(M)}^{2}+O(\varepsilon)^{2}.\end{eqnarray*}
Summarizing, we then expect that\[
\Phi(x_{1},\ldots,x_{n})=\sum_{j=1}^{n}\Vert\xi_{j}\Vert_{L^{2}(M)}^{2}.\]

\subsubsection{Definition of $\Phi^{*}(x_{1},\ldots,x_{n})$.}

This leads us to take the non-rigorous formula for $\Phi$ as a definition
of the non-microstates free Fisher information:

\begin{defn}
\cite{dvv:entropy5} Let $(x_{1},\ldots,x_{n})$ be a family of non-commutative
random variables in $(A,\tau)$. If conjugate variables $(\xi_{1},\ldots,\xi_{n})$
to this family exist, then we set\[
\Phi^{*}(x_{1},\ldots,x_{n})=\sum_{j=1}^{n}\Vert\xi_{j}\Vert_{L^{2}(M)}^{2},\qquad M=W^{*}(x_{1},\ldots,x_{n}).\]
If the conjugate variables do not exist, we set $\Phi^{*}(x_{1},\ldots,x_{n})=+\infty$.
\end{defn}
Note that this definition does not involve microstates.

In the case of a single variable, $\xi_{1}$ ends up being the restriction
of the Hilbert transform of the distribution of $x_{1}$ to the support
of this distribution. One can then compute that if $\mu_{x}$ is Lebesgue
absolutely-continuous, and $d\mu_{x}(t)=p(t)dt$, then\[
\Phi^{*}(x)=\Phi(x)=\frac{2}{3}\int p(t)^{3}dt.\]

\subsubsection{Definition of $\chi^{*}$.}

Since $\Phi^{*}$ was supposed to be proportional to the derivative
of free entropy one can recover free entropy from the free Fisher
information. The formula is\[
\chi^{*}(x_{1},\ldots,x_{n})=\frac{1}{2}\int_{0}^{\infty}\left(\frac{n}{1+t}-\Phi^{*}(x_{1}^{t},\ldots,x_{n}^{t})\right)dt+n\log2\pi e;\]
here as before $x_{j}^{t}=x_{j}+\sqrt{t}s_{j}$, and $(s_{1},\ldots,s_{n})$
is a free semicircular family, free from $(x_{1},\ldots,x_{n})$.

Voiculescu proved that the function\[
t\mapsto\Phi^{*}(x_{1}^{t},\ldots,x_{n}^{t})\]
is monotone decreasing and right semi-continuous in the sense that\[
\lim_{s\to t^{+}}\Phi^{*}(x_{1}^{s},\ldots,x_{n}^{s})=\Phi^{*}(x_{1}^{t},\ldots,x_{n}^{t}).\]
It is an important open question if this function is always continuous. 

Furthermore, if $n=\sum\tau(x_{j}^{2})$, then\[
\frac{n}{1+t}\leq\Phi^{*}(x_{1}^{t},\ldots,x_{n}^{t})\leq\frac{n}{t},\]
which implies that the integral defining $\chi^{*}$ makes sense and
converges to a value in $[-\infty,+\infty)$.

\subsection{Properties of $\chi^{*}$.}

As we mentioned in the foreword to this section, the principal outstanding
question in the theory of free entropy is the question of when $\chi=\chi^{*}$.
To this end there are two results:

\begin{itemize}
\item \cite{dvv:entropy5} In the single-variable case, the two quantities
are equal: $\chi(x_{1})=\chi^{*}(x_{1})$;
\item \cite{guionnet-biane-capitaine:largedeviations} In general, the following
inequality is satisfied:\[
\chi(x_{1},\ldots,x_{n})\leq\chi^{*}(x_{1},\ldots,x_{n}).\]

\end{itemize}
The non-microstates definition turns out to be easier to work with
in some respects, but harder in others. One of the big difficulties
in the non-microstates framework is one's inability to prove the change
of variables formula. This difficulty is related to our inability
to handle the continuity properties of the {}``non-commutative Hilbert
transform'', $(x_{1},\ldots,x_{n})\mapsto(\xi_{1},\ldots,\xi_{n})$,
where $(\xi_{1},\ldots,\xi_{n})$ are the conjugate variables to $(x_{1},\ldots,x_{n})$.

Nonetheless, $\chi^{*}$ has a lot of nice properties, for example:
(all of these are from \cite{dvv:entropy5})

\begin{itemize}
\item $\chi^{*}(x_{1},\ldots,x_{n},y_{1},\ldots,y_{m})=\chi^{*}(x_{1},\ldots,x_{n})+\chi^{*}(y_{1},\ldots,y_{m})$
if the families $(x_{1},\ldots,x_{n})$ and $(y_{1},\ldots,y_{m})$
are free;
\item $\chi^{*}(x_{1},\ldots,x_{n},y_{1},\ldots,y_{m})\leq\chi^{*}(x_{1},\ldots,x_{n})+\chi^{*}(y_{1},\ldots,y_{m})$;
\item $\chi^{*}(x_{1},\ldots,x_{n})$ is maximal subject to $\sum\tau(x_{i}^{2})=n$
if and only if $x_{1},\ldots,x_{n}$ are free semicircular variables,
and $\tau(x_{1}^{2})=\cdots=\tau(x_{n}^{2})=1$.
\item If $s_{1},\ldots,s_{n}$ are free semicircular variables, free from
the family $x_{1},\ldots,x_{n}$, then for any $\varepsilon>0$, $\chi^{*}(x_{1}+\sqrt{\varepsilon}s_{1},\ldots,x_{n}+\sqrt{\varepsilon}s_{n})>-\infty$.
\end{itemize}
Comparing the last property of $\chi^{*}$ with the corresponding
property of $\chi$ explains why $\chi=\chi^{*}$ would imply a positive
answer to Connes' embedability question.

\subsection{Non-microstates free entropy dimension.}

Although we don't know how to formulate the packing number definition
of free entropy dimension in the non-microstates approach, the Minkowski
dimension definition does have a straightforward analog. We set\[
\delta^{*}(x_{1},\ldots,x_{n})=n-\liminf_{\varepsilon\to0}\frac{\chi^{*}(x_{1}^{\varepsilon},\ldots,x_{n}^{\varepsilon})}{\log\varepsilon^{1/2}},\]
where as before $x_{j}^{\varepsilon}=x_{j}+\sqrt{\varepsilon}s_{j}$,
and $s_{1},\ldots,s_{n}$ is a free semicircular family, free from
the family $x_{1},\ldots,x_{n}$.

It is tempting to formally apply L'Hopital's rule in the definition
of $\delta^{*}$ and use the fact that $\frac{d}{d\varepsilon}\chi^{*}(x_{1}^{\varepsilon},\ldots,x_{n}^{\varepsilon})=\Phi^{*}(x_{1},\ldots,x_{n})$.
Thus we write\[
\delta^{\star}=n-\liminf_{\varepsilon\to0}\varepsilon\Phi^{*}(x_{1}^{\varepsilon},\ldots,x_{n}^{\varepsilon}).\]
One can easily show that\[
\delta^{\star}(x_{1},\ldots,x_{n})\geq\delta^{*}(x_{1},\ldots,x_{n}),\]
with no examples in which equality does not hold.

There are unfortunately preciously few computations of $\delta^{*}$
or $\delta^{\star}$, and much less is known about their properties
than about the properties of $\delta$. In particular, it is not known
in general if $\delta^{\star}$ or $\delta^{*}$ depend only on the
algebra generated by $x_{1},\ldots,x_{n}$, taken with its trace.

We summarize what is known below:

\begin{itemize}
\item $\delta^{*}(x_{1})=\delta^{\star}(x_{1})=\delta_{0}(x_{1})=1-\sum_{t}\mu_{x_{1}}(\{ t\})^{2}$,
where $\mu_{x}$ is the law of $x_{1}$;
\item $\delta^{*}(x_{1},\ldots,x_{n},y_{1},\ldots,y_{m})=\delta^{*}(x_{1},\ldots,x_{n})+\delta^{*}(y_{1},\ldots,y_{m})$
if $(x_{1},\ldots,x_{n})$ are free from $(y_{1},\ldots,y_{m})$;
the same is true for $\delta^{\star}$;
\item \cite{connes-shlyakht:l2betti} If $x_{1},\ldots,x_{n}$ are generators
of a tracial algebra $(A,\tau)$, then $\delta^{*}(x_{1},\ldots,x_{n})\leq\delta^{\star}(x_{1},\ldots,x_{n})\leq\beta_{1}^{(2)}(A,\tau)-\beta_{0}^{(2)}(A,\tau)+1,$
where $\beta_{j}^{(2)}(A,\tau)$ are the $L^{2}$-Betti numbers of
$(A,\tau)$.
\item \cite{shlyakht-mineyev:freedim} If $x_{1},\ldots,x_{n}\in\mathbb{C}\Gamma$
are self-adjoint and generate the group algebra of a discrete group
$\Gamma$, then equality holds: $\delta^{\star}(x_{1},\ldots,x_{n})=\delta^{*}(x_{1},\ldots,x_{n})=b_{1}^{(2)}(\Gamma)-b_{0}^{(2)}(\Gamma)+1$,
where $b_{j}^{(2)}(\Gamma)=\beta_{j}^{(2)}(\mathbb{C}\Gamma)$ are
the $L^{2}$-Betti numbers of $\Gamma$. In particular, in this case
$\delta^{*}=\delta^{\star}$ are algebraic invariants;
\item \cite{connes-shlyakht:l2betti} If $W^{*}(x_{1},\ldots,x_{n})$ has
diffuse center, then $\delta^{\star}(x_{1},\ldots,x_{n})\leq1$.
\end{itemize}
\bibliographystyle{amsalpha}
\providecommand{\bysame}{\leavevmode\hbox to3em{\hrulefill}\thinspace}

\end{document}